\documentclass{amsart}

\newtheorem{theorem}{Theorem}[section]

\newtheorem{lemma}[theorem]{Lemma}

\newcommand{\ZZ}{{\mathbb Z}}

\newcommand{\II}{{\mathbb I}}

\newcommand{\cP}{{\mathcal P}}
\newcommand{\cQ}{{\mathcal Q}}

\title[ Limit periodic Jacobi matrices
]{  Limit periodic Jacobi matrices 
  with a prescribed p-adic
hull and a singularly continuous spectrum  }
\author{F. Peherstorfer, A. Volberg, P. Yuditskii}

\begin{document}

\maketitle
\begin{abstract}
In this work we build a certain machine that allows
to construct almost periodic Jacobi matrices with singularly
continuous spectrum 
a prescribed p-adic
hull.
\end{abstract}

\section{Introduction}
A Jacobi matrix $J:l^2(\ZZ)\to l^2(\ZZ)$ is called almost periodic if
the family
$$
\{S^{-k}JS^k\}_{k\in \ZZ},
$$
where $S$ is the shift operator, $S|k\rangle=|k+1\rangle$,
is a precompact in the operator topology.

\begin{proof}[Example]
Let $G$ be a compact abelian group, $p(\alpha), q(\alpha)$ be continuous
functions on $G$,  $p(\alpha)\ge 0$. Then $J(\alpha)$ with the
coefficient sequences $\{p(\alpha+k\mu)\}_k, \{q(\alpha+k\mu)\}_k$,
$\mu\in G$, is almost periodic.
\end{proof}

Let us show that in fact this is a general form of almost periodic
Jacobi matrices. For a given almost periodic $J$ define the metric on
$\ZZ$ by
$$
\rho_J(k):=||S^{-k}JS^k-J||.
$$
Evidently $\rho_J(k+m)\le\rho_J(k)+\rho_J(m)$.  Then $J=J(0)$,
where $G=I_J$, $I_J$ is the closer
of $\ZZ$ with respect to $\rho_J$, and $\mu=1\in I_J$.

Recall that for a given system of integers $\{d_k\}_{k=1}^\infty$
one can define
 the set 
\begin{equation}
\II=\underleftarrow{\lim}\{\ZZ/d_1...d_k\ZZ\},
\end{equation}
 that is $\alpha\in \II$ means that $\alpha$ is
 a sequence $\{\alpha_0,\alpha_1,\alpha_2,...\}$ such that
 $$
 \alpha_k\in \ZZ/d_1...d_{k+1}\ZZ\quad \text{and}
 \quad \alpha_k|\text{mod} \,d_1...d_k=\alpha_{k-1}.
 $$
 In particular, if $p$ is a prime number and
 $d_k=p$ we get the ring of $p$--adic integers, $\II=\ZZ_p$.
 
In this work we build a certain machine that allows
to construct almost periodic Jacobi matrices with singularly
continuous spectrum such that $I_J=\II$.

Let
$T$
be an expanding polynomial, $\deg T= d$.
Under the normalization 
$$
T^{-1}:[-1,1]\to [-1,1]
$$
such a polynomial is well defined by position
of its critical values 
$$
\{t_i=T(c_i): T'(c_i)=0,\  c_i>c_j \ \text{for}\ i>j
\}.
$$
The key element of the construction is the following 

\begin{theorem} 

Let $\tilde J$ be a Jacobi matrix
with the spectrum on $[-1,1]$. Then the following
Renormalization Equation has a solution 
$J=J(\epsilon,\tilde J)=J(\epsilon,\tilde J; T)$ with
the spectrum on $T^{-1}([-1,1])$:
\begin{equation}\label{3}
V^*_{\epsilon}(z-J)^{-1}V_
{\epsilon}=
(T(z)-\tilde J)^{-1}T'(z)/d,
\end{equation}
where $V_\epsilon|k\rangle=|\epsilon+dk\rangle$, 
$0\le\epsilon\le d-1$.
Moreover, 
 if
$\min_i |t_i|\ge 10$ then 
$$
||J(\epsilon,\tilde J_1)-J(\epsilon,\tilde J_2)||\le
c||\tilde J_1-\tilde J_2||.
$$
with an absolute constant $c<1$ (does not depend of $T$ also).
\end{theorem}

Let us point out the following two properties
of the function $J(\epsilon,\tilde J; T)$. First,  due to the
commutant relation  $V_\epsilon S=S^d V_\epsilon$ one gets
${J(\epsilon,S^{-m}\tilde J S^m)=S^{-d m}J(\epsilon,\tilde J)S^{d m}}$.
Second, the chain rule holds
$$
J(\epsilon_0, J(\epsilon_1,\tilde J; T_2); T_1)=
J(\epsilon_0+\epsilon_1 d_1,\tilde J; T_2\circ T_1),
$$
where $d_i=\deg T_i$, $0\le\epsilon_i\le d_{i+1}$.

Next steps are quite simple.
For given $d_1, d_2...$, let us chose polynomials
$T_1, T_2...$, $\deg T_k=d_k$ with sufficiently large
critical values. For a fixed sequence $\epsilon_0, \epsilon_1...$,
$0\le\epsilon_k\le d_{k+1}$, 
define
$J_n=J(\epsilon_0+\epsilon_1 d_1+...+\epsilon_{n-1} d_1...d_{n-1},
\tilde J; T_n\circ ...\circ T_2\circ T_1)$. Then $J=\lim_{n\to\infty}
J_n$ exists and does not depend of $\tilde J$. Moreover,
$$
||J-S^{-d_1...d_l m}JS^{d_1...d_l m}||\le A c^l,
\ A>0.
$$
That is $\rho_J$  defines on $\ZZ$ the standard 
$p$--adic topology  in this
case.

\section{
Renormalization equation} 

Let
  $T(z)=z^d+...$ be
 an expanding polynomial,
 $T^{-1}:[-\xi,\xi]\to [-\xi,\xi]$. If $\tilde J:l^2(\ZZ)\to l^2(\ZZ)$
is a Jacobi matrix with the spectrum on $[-\xi,\xi]$ we define
$J$ as a Jacobi matrix with the spectrum on
$T^{-1}([-\xi,\xi])$ that satisfied the Renormalization Equation
\begin{equation}\label{t01}
V^*(z-J)^{-1}V=(T(z)-\tilde J)^{-1}T'(z)/d,
\quad V|k\rangle=|kd\rangle,
\end{equation}

\begin{lemma}\label{1.1} Let $J^{(s)}$ be the
$s$-th $d\times d$ block of the matrix $J$, that is
\begin{equation}
J^{(s)}=
\begin{bmatrix}
q_{sd}&p_{sd+1}& & & \\
p_{sd+1}&q_{sd+1}&p_{sd+2}& & \\
    &\ddots&\ddots&\ddots & \\
& & p_{sd+d-2}&q_{sd+d-2}&p_{sd+d-1} \\
&  & &p_{sd+d-1}&q_{sd+d-1}
\end{bmatrix}.
\end{equation}
Then its resolvent function is of the form
\begin{equation}\label{3}
\left<0\left|(z-J^{(s)})^{-1} \right|0\right>=
\frac{T'(z)/d}{T^{(s)}(z)}.
\end{equation}
Moreover at the critical points $\{c: T'(c)=0\}$ the following
decomposition in a continued fraction holds true
\begin{equation*}
T^{(s)}(c)=T(c)-\tilde q_s- \frac{\tilde p^2_s} {T(c)-\tilde
q_{s-1}-...},
\end{equation*}
that is
\begin{equation}\label{4}
\frac 1{T^{(s)}(c)}=\langle s|(T(c)-\tilde J_-(s))^{-1}
|s\rangle,
\end{equation}
where $\tilde J_-(s)=P_{l^2_-(s)}\tilde J|l^2_-(s)$,
$l^2_-(s)$ is formed by $\{|s+k\rangle\}_{k\le 0}.$ Also
\begin{equation}\label{4.1}
\tilde p_{s+1}=p_{sd+1}...p_{sd+d}
\quad\text{and}\quad q_{sd}=\tilde q_s.
\end{equation}
\end{lemma}
Note that
\begin{equation}\label{4.2}
T^{(s)}(z)=(z-\tilde q_s)T'(z)/d+
\sum_{c:T'(c)=0}\frac{T'(z)}{(z-c)T''(c)}T^{(s)}(c)
\end{equation}
that is we can restore $J^{(s)}$ by \eqref{4},
\eqref{3} and then find $p_{sd+d}$ by \eqref{4.1}.
\begin{proof} 
We write $J$ as a $d\times d$ block matrix (each block is of
infinite size):
\begin{equation}\label{5}
J=
\begin{bmatrix}
\cQ_{0}&\cP_{1}& & &S\cP_d  \\
\cP_{1}&\cQ_{1}&\cP_{2}& & \\
    &\ddots&\ddots&\ddots & \\
& & \cP_{d-2}&\cQ_{d-2}&\cP_{d-1} \\
\cP_d S^*&  & &\cP_{d-1}&\cQ_{d-1}
\end{bmatrix}.
\end{equation}
Here $\cP_k$ (respectively $\cQ_k$) is a diagonal matrix
$\cP_k=\text{diaq}\{p_{k+sd}\}_{s\ge 0}$ and $S$ is the
shift operator.  In this case $V^*$ is the projection on the
first block--component.

Using this representation and being well
    known identity for  block matrices
\begin{equation*}
\begin{bmatrix}
A&B\\
C&D
\end{bmatrix}^{-1}=
\begin{bmatrix}
(A-B D^{-1}C)^{-1}&*\\
*&*
\end{bmatrix},
\end{equation*}
we get
\begin{equation}\label{6}
\frac{T(z)-\tilde J}{T'(z)/d}= z-\cQ_0-
\begin{bmatrix}\cP_1,&...,&S\cP_d
\end{bmatrix}
\{z-J_{1}\}^{-1}
\begin{bmatrix}
\cP_1\\
\vdots\\
\cP_d S^*
\end{bmatrix},
\end{equation}
where $J_{1}$ is the matrix that we obtain from $J$ by deleting
the first block--row and the first block--column in \eqref{5}.
Note that in $(z-J_{1})$ each block is a diagonal matrix that's
why we can easily get an inverse matrix in terms of orthogonal
polynomials.

Let us introduce the following notations: everything related to
$J^{(s)}$ has superscript $s$. For instance: $p_k^{(s)}=p_{sd+k}$,
$1\le k\le d$, respectively $P^{(s)}_d$ and $Q^{(s)}_d$ mean
orthonormal polynomials of the first and second kind. In this
terms equation \eqref{6} is equivalent to the two series of scalar
relations corresponding to the diagonal and off diagonal entries
\begin{equation}\label{7}
\frac{T(z)-\tilde q_{s+1}}{T'(z)/d}= \frac{P_d^{(s+1)}(z)}
{Q^{(s+1)}_d(z)}- p_{ds}^2\frac{Q_{d-1}^{(s)}(z)/p_{ds}}
{Q^{(s)}_d(z)}
\end{equation}
     and
\begin{equation}\label{8}
\frac{\tilde p_{s+1}}{T'(z)/d}= \frac{p_1^{(s)} ... p_d^{(s)}}
{z^{d-1}+...}=\frac 1{Q^{(s)}_d(z)}.
\end{equation}
We have to remind (see \eqref {3} and \eqref{8}) that
\begin{equation*}
\frac{Q_d^{(s)}(z)} {P^{(s)}_d(z)} =\frac{z^{d-1}+...}{z^d+...} =
\frac{T'(z)/d} {T^{(s)}(z)}.
\end{equation*}

    Now, due to the
Wronskian identity, if $T'(c)=0$ then
\begin{equation}\label{9}
-p_{ds}{Q^{(s)}_{d-1}(c)} =\frac 1{P^{(s)}_d(c)}.
\end{equation}
So, combining \eqref{7}, \eqref{8} and \eqref{9} we get the
recurrence relation
\begin{equation}
T(c)-\tilde q_{s+1}= T^{(s+1)}(c)+\frac{\tilde p^2_{s+1}}
{T^{(s)}(c)}.
\end{equation}

\end{proof}

Let us mention that the Renormalization Equation can be rewritten
equivalently in the  form of polynomials equations.

\begin{lemma} Equation \eqref{t01} is equivalent to
\begin{equation}
\label{re.1}
V^*T(J)=\tilde JV^*,
\end{equation}
\begin{equation}
\label{re.2}
V^*\frac{T(z)-T(J)}{z-J}V=T'(z)/d.
\end{equation}
\end{lemma}

\begin{proof}
Starting with \eqref{re.1}, \eqref{re.2} we get
\begin{equation*}
(T(z)-\tilde J)V^*(z-J)^{-1}V=
V^*\{T(z)-T(J)\}\{z-J\}^{-1}V=T'(z)/d.
\end{equation*}

Having \eqref{t01} we get
\begin{equation}\label{re.p.1}
\begin{split}
V^*\frac{T(z)-T(J)}{z-J}V=&
T(z)V^*(z-J)^{-1}V-V^*\frac{T(J)}{z-J}V\\
=&T(z)\frac{T'(z)/d}{T(z)-\tilde J}-V^*\frac{T(J)}{z-J}V\\
=&{T'(z)/d}+\tilde JV^*(z-J)^{-1}V-V^*\frac{T(J)}{z-J}V.
\end{split}
\end{equation}
Since the left hand side in \eqref{re.p.1} is a polynomial of $z$
we obtain two relations
\begin{equation*}
V^*\frac{T(z)-T(J)}{z-J}V=T'(z)/d
\end{equation*}
and
\begin{equation*}
\{\tilde JV^*-V^*T(J)\}\{(z-J)^{-1}V\}=0.
\end{equation*}
Since vectors of the form $(z-J)^{-1}V f$, $f\in l^2$, are compleat
in $l^2$ the last relation implies \eqref{re.1}.
\end{proof}

\section{Proof of the main theorem}

We start with   (undoubtedly well--known and  simple)
\begin{lemma}\label{3.1}
Assume that two non--normalized measures $\sigma$ and $\tilde\sigma$
are mutually absolutely continuous. Moreover,
 $d\tilde\sigma
=f\,d\sigma$ and 
$(1+\epsilon)^{-1}\le f\le (1+\epsilon)$.
    Let us associate with these measures Jacobi matrices
$J= J(\sigma)$,
$\tilde J= J(\tilde\sigma)$.
 Then for their coefficients we have
$$
|\tilde p_s-p_s|\le \epsilon|| J||, \quad s\ge 0.
$$
\end{lemma}

\begin{proof} Assume that $p_s\ge\tilde p_s$.
Let us use an extreme property of orthogonal polynomials,
\begin{equation*}
\begin{split}
(1+\epsilon)\tilde p_0^2...\tilde p_s^2=& (1+\epsilon)\int\tilde
p_0^2...\tilde p_s^2
\tilde P_s^2\,d\tilde\sigma \ge
\int\{z^s+...\}^2\,d \sigma\\
\ge& \inf_{\{P=z^s+...\}}\int P^2 \,d \sigma
= p_0^2... p_s^2.
\end{split}
\end{equation*}
Similarly
\begin{equation*}
(1+\epsilon)p_0^2... p_{s-1}^2
    \ge
\tilde p_0^2...\tilde p_{s-1}^2.
\end{equation*}
Therefore
$$
 p_s^2\le \tilde p_s^2\le{(1+\epsilon)^2}p_s^2
$$
    and hence
$$
0\le \tilde p_s-p_s\le \epsilon p_s.
$$

\end{proof}

\begin{proof}[Proof of the Theorem]
Given $\tilde J_1$ and $\tilde J_2$ let us compere the blocks
$J^{(s)}_1$ and $J^{(s)}_2$ of the matrix $J_1$ and $J_2$.
To this end consider the function
$$
f(c):=\frac{T^{(s)}_2(c)}{T^{(s)}_1(c)}
$$
Assuming $f(c)\ge 1$ let us estimate $f(c)-1$ from above.

\begin{equation}
\begin{split}
f(c)-1=&\frac{1/T^{(s)}_1(c)-1/T^{(s)}_2(c)}{1/T^{(s)}_2(c)}\\
=&
\frac{\langle s|(T(c)-\tilde J_{2,-})^{-1}
(\tilde J_{1,-}-\tilde J_{2,-})(T(c)-\tilde
J_{1,-})^{-1}|s\rangle}{1/T^{(s)}_2(c)}.
\end{split}
\end{equation}
Since
$$
|1/{T^{(s)}_2(c)}|=\left|\int_{-\xi}^\xi
\frac{d\sigma(x)}{T(c)-x}
\right|\ge\frac 1{|T(c)|+\xi},
$$
and
$$
||
(T(c)-\tilde J_{i,-})^{-1}|s\rangle
||\le\frac 1{|T(c)|-\xi}
$$
we get
\begin{equation}\label{ef}
0\le f(c)-1\le ||\tilde J_1-\tilde J_2||\frac
{|T(c)|+\xi}{(|T(c)|-\xi)^2}.
\end{equation}
Thus, by Lemma \ref{3.1}, we obtain
\begin{equation}\label{e1}
|(p_1)_{sd+k}-(p_2)_{sd+k}|\le\delta ||\tilde J_1-\tilde J_2||,
\quad 1\le k\le d-1,
\end{equation}
where
$$
\delta:=\max_{c}\frac
{|T(c)|/\xi+1}{(|T(c)|/\xi-1)^2}.
$$

We have to estimate $|(p_1)_{sd+d}-(p_2)_{sd+d}|$.
First we claim that
\begin{equation}\label{14}
\frac 1{(p_i)_{sd+1}...(p_i)_{sd+d-1}}\le
\max_c\frac 1{|T(c)|/\xi-1}.
\end{equation}
Recall that
$$
\frac{1}{z-q_{sd+1}+...}=\frac {A(z)}{B(z)}=
\int\frac{d\sigma(x)}{z-x},\quad\int\,d\sigma=1,
$$
where (see \eqref{3}, \eqref{4.2})
$$
B(z)=\frac{T'(z)/d}{p_{sd+2}...p_{sd+d-1}},\quad A(z)=\frac
1{p_{sd+1}^2}
\frac{T^{(s)}(z)-(z-q_s)T'(z)/d}{ p_{sd+2}...p_{sd+d-1}}.
$$
Let $B_1$, $A_1$ be orthonormal polynomials of the first and second
kind, $\deg B_1=d-2$, so that
$$
B_1(z) A(z)-A_1(z)B(z)=1,
\quad
\int B_1^2 d\sigma=1.
$$
Since at the critical points $B_1(c)=1/A(c)$ we get
$$
\int\frac{{p_{sd+1}^2}}{T^{(s)}(c)^2}d\sigma=
\frac 1{ p^2_{sd+1}...p^2_{sd+d-1}}.
$$
Thus \eqref{14} is proved.

Now, by \eqref{4.1}
\begin{equation*}
\begin{split}
(p_1)_{sd+d}-(p_2)_{sd+d}=&
\frac{(\tilde p_1)_{s+1}}{ (p_1)_{sd+1}...(p_1)_{sd+d-1}}-
\frac{(\tilde p_2)_{s+1}}{ (p_2)_{sd+1}...(p_2)_{sd+d-1}}\\
=&
\frac{(\tilde p_1)_{s+1}-
(\tilde p_2)_{s+1}}{ (p_1)_{sd+1}...(p_1)_{sd+d-1}}\\
+&\frac{(\tilde p_2)_{s+1}}{ (p_1)_{sd+1}...(p_1)_{sd+d-1}}
\left(1-
\frac{ (p_1)_{sd+1}...(p_1)_{sd+d-1}}
{ (p_2)_{sd+1}...(p_2)_{sd+d-1}}
\right)
\end{split}
\end{equation*}
Using \eqref{14}, \eqref{ef} and Lemma \ref{3.1}
we obtain
\begin{equation}\label{e2}
\begin{split}
|(p_1)_{sd+d}-(p_2)_{sd+d}|\le&
||\tilde J_1-\tilde J_2||\max_c\frac 1{|T(c)|/\xi-1}\\
+&\max_c\frac 1{|T(c)|/\xi-1}||\tilde J_1-\tilde J_2||\delta/2\\
=&||\tilde J_1-\tilde J_2||\max_c\frac 1{|T(c)|/\xi-1}
(1+\delta/2).
\end{split}
\end{equation}
Thus \eqref{e1} and \eqref{e2} show that say for
$\min_c{|T(c)|/\xi}\ge 10$ the renormalization is a contraction.
\end{proof}
\end{document}